\documentclass{article}
\usepackage{arxiv}
\usepackage{amsmath,amssymb,amsfonts}
\usepackage{amsthm}
\theoremstyle{definition}
\newtheorem{definition}{Definition}[section]
\newtheorem{proposition}[definition]{Proposition}
\newtheorem{theorem}[definition]{Theorem}
\newtheorem{corollary}[definition]{Corollary}
\newtheorem{remark}[definition]{Remark}

\newtheorem{example}[definition]{Example}

\usepackage[utf8]{inputenc} % allow utf-8 input
\usepackage[T1]{fontenc} % use 8-bit T1 fonts
\usepackage{hyperref} % hyperlinks
\hypersetup{colorlinks}
\usepackage{url} % simple URL typesetting
\usepackage{booktabs} % professional-quality tables
\usepackage{amsfonts} % blackboard math symbols
\usepackage{nicefrac} % compact symbols for 1/2, etc.
\usepackage{microtype} % microtypography
\usepackage{cleveref} % smart cross-referencing
\usepackage{lipsum} % Can be removed after putting your text content
\usepackage{graphicx}
\usepackage{doi}

\title{Ideal and weak Amenability of Fr\'{e}chet locally $C^*$-algebra}

\author{ \href{https://orcid.org/0000-0000-0000-0000}{\includegraphics[scale=0.06]{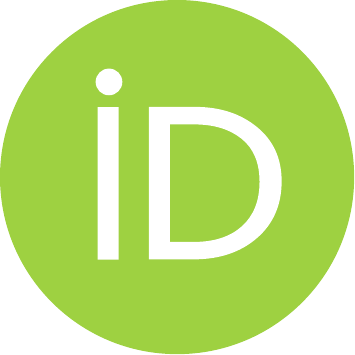}\hspace{1mm}Ali Ranjbari}%
\thanks{A. Ranjbari: \ Tel.: +98-9124528763, \quad Fax: +98-3137934631} \\
	Department of Mathematics\\
	University of Isfahan\\
	 Isfahan, 	 Iran\\
	\texttt{aranjbari1353@yahoo.com} \\
	\And
	\href{https://orcid.org/0000-0000-0000-0000}{\includegraphics[scale=0.06]{orcid.pdf}\hspace{1mm}Ali Rejali} \\
	Department of Mathematics\\
	University of Isfahan\\
	Isfahan, Iran\\
	\texttt{rejali@sci.ui.ac.ir} \\
}

\renewcommand{\headeright}{\thepage}
\renewcommand{\shorttitle}{$n$-ideal and $n$-weak amenability of Fr{\'e}chet algebras}

\hypersetup{
pdftitle={$n$-ideal and $n$-weak amenability of Fr{\'e}chet algebras},
pdfsubject={n-ideal and n-weak amenability of Frechet algebras},
pdfauthor={Ali Ranjbari, Ali Rejali},
pdfkeywords={Frechet algebra, n-ideal amenability, n-weak amenability},
}

\DeclareMathOperator{\ad}{ad}
\DeclareMathOperator{\supp}{supp}

\begin{document}
\maketitle

\begin{abstract}
The notion of Fr\'{e}chet locally $C^*$-algebra generalizes the notion of $C^*$-algebra. 
In this paper, we first present some definitions and basic facts about locally $C^*$-algebra, and then we introduce and study the notion of ideal a and weak menability for these algebras. Also,we show that every Fr\'{e}chet locally $C^*$-algebra is ideally amenable.
\end{abstract}

% keywords can be removed
\keywords{Fr\'{e}chet Algebras, Inverse limits of $C^*$-algebras, Locally convex algebra, Ideal amenability of Fr\'{e}chet locally $C^*$-algebra.}
\textbf{SubClass}: 46H05, 46J05, 46A04, 43A15, 43A60.

\section{Introduction}
Some of the notions related to Banach algebras have been introduced and studied for Fr\'{e}chet algebras. For example, in \cite{13}, the contractibility of Fr\'{e}chet algebras was investigated. 
Moreover, the notion of amenability of a Fr\'{e}chet algebra was introduced by Helemskii \cite{14} and studied by Pirkovskii \cite{19}. 
Also in \cite{17}, Lawson and Read introduced and studied some notions about approximate amenability and approximate contractibility of Fr\'{e}chet algebras. Rejali and et al in \cite{1}, introduced and studied the notion of weak amenability of Fr\'{e}chet algebra. Moreover,they introduced the notion of Segal Fr\'{e}chet algebra in \cite{3} and the semisimple Fr\'{e}chet algebra in \cite{2}, recently, they introduce the concepts of $\phi$-amenability and character amenability for Fr\'{e}chet algebras in \cite{4}.
Furthermore, in \cite{22} we introduced and studied the notion of Lipschitz algebra for Fr\'{e}chet algebra
\\
Moreover, we introduced the notion of I-weak amenability and also ideal amenability for Fr\'{e}chet algebra in \cite{21}. recently, in \cite{23} we introduced the notion of n-ideal amenability and n-weak amenability for Fr\'{e}chet algebra and we examine how these concepts in Banach algebra can be generalized and defined for Fr\'{e}chet algebra. Pirkovskii \cite{19} extend a number of characterizations of amenability obtained by Johnson and by Helemskii to the setting of locally $m$-convex Fr\'{e}chet algebras. As a corollary, he shows that Connes and Haagerup's theorem on amenable $\mathbf{C^*}$-algebra and Sheinberg's theorem on amenable uniform algebras hold in the Fr\'{e}chet algebra case. 

\noindent The notion of Fr\'{e}chet locally $\mathbf{C^*}$-algebra generalizes the notion of $\mathbf{C^*}$-algebra. In this paper, we first present some definitions and basic facts about locally $\mathbf{C^*}$-algebra, and then we introduce and study the notion of ideal amenability for Fr\'{e}chet locally $C^*$-algebra.

\section{preliminaries}
We first recall some of the basic facts about Fr\'{e}chet algebra, some of their properties and locally $\mathbf{C^*}$-algebra; for further details, see \cite{7}, \cite{18} and \cite{24}.

\noindent A locally convex space $E$ is a topological vector space in which each point has a neighborhood basis of convex sets. Throughout the paper, all locally convex spaces are assumed to be Hausdorf.

\noindent A Fr\'{e}chet space is a metrizable, complete and locally convex vector space.

\noindent A Fr\'{e}chet algebra is a complete algebra $A$ generated by a sequence $\left(p_n\right)_{n\in \mathbb{N}}$ of separating increasing submultiplicative seminorms, i.e 
\\
$p_n(xy)\le p_n(x) p_n(y)$ for all $n \in \mathbb{N}$ and every $x; y \in A$ such that $p_n(x)\le \ p_{n+1}(x)$ for all positive integer $n$ and $x\in A$. 
We write $\left(A,(p_n)\right)$ for the corresponding Fr\'{e}chet algebra.
Closed subspaces of Fr\'{e}chet spaces are Fr\'{e}chet and quotients of Fr\'{e}chet spaces by closed subspaces are Fr\'{e}chet. Limits of countable many Fr\'{e}chet spaces are Fr\'{e}chet. The Fr\'{e}chet spaces are exactly the (projective) limits of sequences of Banach spaces.
\\
An involution on an algebra $A$ is a map $*\colon A \to A$, $x \mapsto x^*$, such that for all $x$, $y\in A$ and $\lambda \in \mathbb{C}$, 
\begin{enumerate}
\item[(i)]
$(x + y)^* = x^*+ y^*$,
\item[(ii)]
$(\lambda x)^* = \lambda x^*$,
\item[(iii)]
$(xy)^* = y^*x^*$,
\item[(iv)]
$(x^*)^*=x$.
\end{enumerate}
A topological algebra A equipped with an involution $*$, this will always be assumed to be continuous, is called a topological $*$-algebra. 
The symbol $A[\tau]$ will denote a topological Because, whose given topology is $\tau$.

Let $A$ be a $*$-algebra. A seminorm $p$ on $A$ is called $m^*$-seminorm, resp.

$C^*$-seminorm, if it is submultiplicative, i.e $p\left(xy\right)\leq p\left(x\right) p\left(y\right)$ for all $x; y \in A$ and $p(x^*)=p(x)$, resp. $p(xx^*) = p(x)^2$.

A locally $C^*$-algebra is a complete Hausdorff complex $*$-algebra $A$ whose topology is determined by its continuous $\boldsymbol{C^*}$-seminorms in the sense that a net $(a_\alpha)_{\alpha \epsilon \Lambda}$ converges to $0$ if and only if the net $p(a_\alpha)_{\alpha \epsilon \Lambda}$ converges to $0$ for each continuous $\boldsymbol{C^*}$-seminorm $p$ on $A$. 
A Fr\'{e}chet locally $\boldsymbol{C^*}$-algebra is a locally $C^*$-algebra whose topology is determined by a countable family of $\boldsymbol{C^*}$-seminorms.
Clearly, any $\boldsymbol{C^*}$-algebra is a Fr\'{e}chet locally $\boldsymbol{C^*}$-algebra.

Nuclear $\boldsymbol{C^*}$-algebra is a $\boldsymbol{C^*}$-algebra $A$ such that the injective and projective $C^*$-cross norms on $A\otimes B$ are the same for every $\boldsymbol{C^*}$-algebra $B$. 
This property was first studied by Takesaki (1964) under the name "Property T", which is not related to Kazhdan's property T. 

Let $\left(A,(p_k)\right)$ be a Fr\'{e}chet algebra and $X$ be a locally convex $A$-bimodule. 
A continuous derivation of $A$ into $X$ is a continuous linear mapping $D\colon A \to X$ such that for all $a, b \in A$, $D(ab) = a \cdot D(b) + D(a) \cdot b$. 
If $x \in X$, we define $\ad_x\colon A \to X$ by $\ad_x(a) = a \cdot x - x\cdot a$ for each $a \in A$. 
Then $\ad_x$ is a derivation. Such derivations are called inner derivations. The Fr\'{e}chet algebra $A$ is amenable if every continuous derivation $D\colon A\to X*$ is inner for every locally convex $A$-bimodule $X$, in the other words $H^1(A, X^*) = 0$ and $A$ is weakly amenable if $H^1(A, A^*) = 0$.

Similar to the Banach algebra case, for every closed (two-sided) ideal $I$ in $A$ we say that $A$ is $I$-weakly amenable if every continuous derivation $D\colon A\to I^*$ is inner. 
Moreover $A$ is called ideally amenable if it is $I$-weakly amenable, for every closed (two-sided) ideal $I$ in $A$.

\section{Ideal Amenability of Fr\'{e}chet locally $\boldsymbol{C^*}$-algebra}
The theory of $C^*$-algebra was founded by Gel'fand and Naimark in 1943, and Subsequently, many of their properties, including types of amenability, were studied by various mathematicians. In \cite{12}, Haagerup proved that, a $C^*$-algebra is amenable if and only if it is nuclear. However, every $C^*$-algebra is weakly amenable \cite{24} and too is ideally amenable \cite{9}. The notion of Fr\'{e}chet locally $C^*$-algebra generalizes the notion of $C^*$-algebra. 
Locally $C^*$-algebra were systemically studied by Inove \cite{24}, N. C. Phillips (under the name of pro- $C^*$-algebra \cite{20}), M. Fragoulopoulo \cite{7} and other mathematicians.
\begin{remark}\label{t.rem.1}
Let $A$ be a locally $C^*$-algebra, $S(A)$ the set of all continuous $C^*$-seminorms on $A$ and $p \in S(A)$. 
We set $N_p = \{a \in A \colon p(a) = 0\}$, then $A_p = A/ N_p$ is a $C^*$-algebra in the norm induced by $p$, i.e. $p (a + N_p) = p (a)$, $a \in A$ (\cite[Theorem 10.24]{7},). 
For $p, q \in S(A)$ with $p \ge q$, the surjective Morphisms $\pi_{pq} \colon A_p \to A_q$ defined by $\pi_{pq} (a+N_p) = a+N_q$ induce the inverse system $\{A_p\; \pi_{pq}\}_{p,q \in S(A), p\ge q}$ of $C^*$-algebras and $A = \varprojlim A_p$, i.e., The Arens-Michael decomposition gives us the representation of $A$ as an inverse limit of $C^*$-algebra.
\end{remark}
\begin{theorem}\label{t.the.1}
A Fr\'{e}chet locally $C^*$-algebra is amenable if and only if it is nuclear. 
\end{theorem}
\begin{proof}
Let $A$ be a nuclear Fr\'{e}chet locally $C^*$-algebra, and let $A = \varprojlim A_p$, be its Arens-Michael decomposition. 
By \cite[4.1]{20} for each continuous $C^*$-seminorm $p$ on $A$, the $C^*$- algebra $A_p$ is nuclear.
So from Connes-Haagerup Theorem (see \cite[6.5.12]{24}), $A_p$ is amenable.
Then by \cite[Theorem 9.7]{7} A is amenable.

Conversely, let $A =\varprojlim A_p$ be amenable. 
Then by \cite[Theorem 9.7]{7}, for each continuous $C^*$-seminorm $p$ the $C^*$-algebra $A_p$ is amenable. 
So from Connes-Haagerup Theorem (see \cite[6.5.12]{24}), $A_p$ is nuclear, hence By \cite[4.1]{7} the Fr\'{e}chet $C^*$-algebra $A$ is nuclear.
\end{proof}
\begin{definition}\label{t.def.1}
 Let $A[\tau]$ be a locally $C^*$-algebra, and $S(A)$ the set of all continuous $C^*$-seminorms on $A$.
 The $*$-algebra $A_b:= a \in A\colon \supp_\alpha(a) <\infty,p_\alpha \in S(A)\}$ is called the bounded part of $A[\tau]$.

By \cite[Theorem 10.23]{7}, $A_b$ is a $C^*$-algebra under the $C^*$-norm $\| a\|_b:=\supp_\alpha (a)$, $a\in A_b$ and it is dense in $A[\tau]$ (\cite[Theorem 10.24]{5}).
\end{definition}
\begin{proposition}
Let $A$ be a Fr\'{e}chet locally $C^*$-algebra.
Then $A$ is amenable if and only if the $C^*$-algebras $A_b$ is amenable. 
\end{proposition}
\begin{proof}
Let $A$ be a Fr\'{e}chet locally $C^*$-algebra.
By \cite[Theorem 4.5]{5} $A$ is nuclear if and only if the $C^*$-algebras ${A}_{b}$ is nuclear . Therefore if $A$ is amenable, from Theorem \ref{t.the.1}. 
$A$ is nuclear. So the $C^*$-algebras ${A}_{b}$ is nuclear . From Connes-Haagerup Theorem, $A_b$ is amenable.

Conversely, if the $C^*$-algebras $A_b$ is amenable, Then $A_b$ is nuclear. 
By \cite[Theorem 4.5]{5} the Fr\'{e}chet locally $C^*$-algebra $A$ is nuclear, therefore from Theorem \ref{t.the.1}, $A$ is amenable.
\end{proof}
\indent
Recall that a locally convex space $E$ is quasinormable if for each $0$-neighborhood $U\subset E$ there exists a $0$-neighborhood $V \subset U$ such that for each $\epsilon> 0$ there exists a bounded set $B \subset E$ satisfying $V \subset B + \epsilon U$. 
Note that any normed space is quasi-normable. In fact we may take $V = B := U$.
\begin{theorem}\label{t.the.2}\label{t.the.3.5}
Let $A$ be a Fr\'{e}chet algebra and $I$ be a quasinormable closed ideal of $A$.
If $I$ and $\frac{A}{I}$ are ideally amenable, then $A$ is also ideally amenable.
\end{theorem}
\begin{proof}
See \cite[Theorem 4.7]{21} for details.
\end{proof}
\begin{theorem}\label{t.the.3.6}\label{t.the.3}
every Fr\'{e}chet locally $\boldsymbol{C^*}$-algebra is ideally amenable.
\end{theorem}
\begin{proof}
Let $A$ be a Fr\'{e}chet locally $C^*$-algebra, and let $A=\varprojlim A_p$, be its Arens-Michael decomposition.
Since the $A_p= A/N_p$ is a Banach $C^*$-algebra By \cite[4.1]{9}, $A_p$ is ideally amenable. 
On the other hand $N_p = \{a \in A\colon p(a) = 0$, is a closed $*$-ideal of $ A_b$ (bounded part of $A$), therefore $N_p$ is a Banach $C^*$-algebra.
So By \cite[4.1]{9} $N_p$ be an ideally amenable and quasinormable closed ideal of $A$. 
Thus Theorem \ref{t.the.3.5} implying that $A$ is ideally amenable.
\end{proof}
\begin{corollary}\label{t.cor.3.7}
Every Fr\'{e}chet locally $\boldsymbol{C^*}$-algebra is weakly amenable.
\end{corollary}
\indent 
Now we provide some remarkable examples related to all of our results of this work.
\begin{example}\label{t.exa.3.8}
Let $A$ be a non-nuclear Fr\'{e}chet locally $\boldsymbol{C^*}$-algebra. 
Then $A$ is ideally amenable which is not amenable.
\end{example}
\begin{example}\label{t.exa.3.9}
The algebra $c^\infty ([a, b])$ of infinitely differentiable complex-valued functions on a (finite) interval $[a, b]$ in $\mathcal{R}$ is not a Banach space.
However, With seminorms
\[ \mu_k(f) = \sup_{0\le i\le k} \sup_{x\in [a,b]} \left| f^{(i)} (x) \right| \] 
and metric
\[ d\left(f,g\right)=\sum^\infty_{k=0} 2^{-k} \frac{\mu_k(f-g)}{\mu_k\left(f-g\right)+1} \] 
$c^\infty([a,b])$ is commutative Fr\'{e}chet algebra. 
From \cite[Example 3.1]{21} it follows that $C^\infty[0,1]$ cannot be a locally amenable.
Therefore By Theorem \ref{t.the.3} $c^\infty([a,b])$ is not locally $C^*$-algebra.
\end{example}

\section*{Acknowledgements}
{\rm The authors would like to thank the Banach algebra
center of Excellence for Mathematics, University of Isfahan.}

%\bibliographystyle{unsrtnat}
%\bibliography{references} %%% Uncomment this line and comment out the ``thebibliography'' section below to use the external .bib file (using bibtex) .


\begin{thebibliography}{1}
\bibitem{1}
F. Abtahi, S. Rahnama and A. Rejali, Semisimple Segal Fr\'{e}chet algebras, Period. Math. Hungar. 71, no. 2, 146--154 (2015)
\bibitem{2}
F. Abtahi, S. Rahnama and A. Rejali, Weak amenability of Fr\'{e}chet algebras, Politehn. Univ. Bucharest Sci. Bull. Ser. A Appl. Math. Phys. 77, no. 4, 93--104 (2015)
\bibitem{3}
F. Abtahi, S. Rahnama and A. Rejali, Segal Fr\'{e}chet algebras, Anal. Math., (2017) DOI 10.1007/s10476-017-0601-y
\bibitem{4}
Abtahi, F., Rahnama, S. and Rejali, A., $\phi$-amenability and character amenability of Fr\'{e}chet algebras, Forum Math., (2018), aop 
\bibitem{5}
S.J. Bhatt and D.J. Karia, Complete positivity, tensor products and $C^*$-nuclearity for inverse limits of $C^*$-algebras, Proc. Indian Acad. Sci. Math. Sci., 101, no. 3, 149--167 (1991) 
\bibitem{6}
H. G. Dales, Banach algebras and automatic continuity, London Math. Soc, Monographs, 24, Oxford: Clarenden Press (2000)
\bibitem{7}
M. Fragoulopoulou, Topological Algebras with Involution, North-Holland Math. Stud. 200, Elsevier Science, Amsterdam (2005)
\bibitem{8}
H. Goldmann, Uniform Fr\'{e}chet Algebras, North-Holland Math. Stud.162, North-Holland, Amsterdam (1990) 
\bibitem{9}
M. E. Gordji and T. Yazdanpanah, Derivations into Duals of Ideals of Banach Algebras, Proc. Indian Acad. Sci., 114, no. 4, 339--403 (2004)
\bibitem{10}
N. Gronbaek, Weak and syclic amenability for non-commutative Banach algebras, Proc. Edin-burg Math. Soc., 35, no. 2, 315--328 (1992)
\bibitem{11}
S. L. Gulick, The bidual of a locally multiplicatively-convex algebra, Pacific J. Math., 17, 71--96 (1966)
\bibitem{12}
U. Haagerup, All nuclear $C^*$-algebra are amenable, Invent. Math., 74, 305--319 (1983) 
\bibitem{13}
R. El Harti, Contractible Fr\'{e}chet algebras, Proc. Amer. Math. Soc., 132, no. 5, 1251--1255 (2004)
\bibitem{14}
A. Y. Helemskii, The Homology of Banach and Topological Algebras, Math. Appl. (Soviet Series), 41, Kluwer Academic, Dordrecht, (1989)
\bibitem{15}
A. Y. Helemskii, Banach and Locally Convex Algebras, Oxford Sci. Publ., Oxford, New York, (1993)
\bibitem{16}
A. Inoue, Locally $C^*$-algebras, Mem. Faculty Sci. Kyushu Univ. (Ser. A), 25, MR 46: 4219, 197--235 (1971)
\bibitem{17}
P. Lawson and C. J. Read, Approximate amenability of Fr\'{e}chet algebras, Math. Proc. Camb. Phil. Soc., 145, 403--418 (2008)
\bibitem{18}
R. Meise and D. Vogt, Introduction to functional analysis, Oxford Science Publications, (1997) 
\bibitem{19}
A. Y. Pirkovskii, Flat cyclic Fr\'{e}chet modules, amenable Fr\'{e}chet algebras, and approximate identities, Homology Homotopy Appl., 11, no. 1, 81--114 (2009)
\bibitem{20}
N. C. Phillips, Inverse limit of $C^*$ -algebra, J. Operator Theory, 19, 159--195 (1988)
\bibitem{21}
Ranjbari, A. and Rejali, A., Ideal amenabilit y of Fr\'{e}chet algebra. U.P.B. Sci. Bull., Series A, 79, Iss. 4, 51--60 (2017)
\bibitem{22}
Ranjbari, A. and Rejali, A., Fr\'{e}chet a-lipschitz vector-valued operator algebras, U.P.B. Sci. Bull., Series A, 80, Iss. 4, (2018)
\bibitem{23}
Ranjbari, A. and Rejali, A., $n$-ideal and $n$-weak amenability of Fr\'{e}chet algebra, Acta Mathematica Sinica, English Series. to Appear.
\bibitem{24}
V. Runde, Lectures on Amenability, Lecture Notes in Mathematics1774, Springer, (2002)




\end{thebibliography}

\end{document}